\documentclass[11pt,twoside,a4paper]{amsart}
\usepackage[english]{babel}

\usepackage[utf8]{inputenc}

\usepackage{pstricks-add}
\usepackage{array}
\usepackage{amsfonts}
\usepackage{amssymb}
\usepackage{amsmath}
\usepackage{graphicx}
\usepackage{graphics}
\usepackage{amsthm}

\usepackage{amssymb}

\usepackage{mathabx}

\usepackage[cmtip,all]{xy}

\usepackage{fancyhdr}
 
\fancypagestyle{noheadrule}{
\fancyhf{}

\fancyhead[LE,RO]{\thepage}
}

\newtheorem{theorem}{Theorem}[section] 

\newtheorem{proposition}[theorem]{Proposition}
\newtheorem{remark}[theorem]{Remark}

\newcommand{\xg}{$(X, \mathbb{G})$ }
\newcommand{\yg}{$(Y, \mathbb{G})$ }
\newcommand{\la}{\mathsf{\Lambda} }
\newcommand{\xlag}{$(X_{\mathsf{\Lambda}}, \mathbb{G})$ }

\normalsize

\thanks{
\\
This work was funded by the FONDECYT Post-doctoral grant n. 3150535.}
\begin{document}
\title[Model sets and almost automorphy]{A short proof of the characterization of model sets by almost automorphy}
\normalsize
\author{Jean-Baptiste Aujogue}

\begin{abstract}
The aim of this note is to provide a conceptually simple demonstration of the fact that repetitive model sets are characterized as the repetitive Meyer sets with an almost automorphic associated dynamical system. 
\end{abstract}

\maketitle

\section{Introduction} 
\label{intro}

In this note we will be interested in certain point sets of a locally compact Abelian group $\mathbb{G}$ that are nowadays called \textit{Meyer sets}. Meyer introduced such point sets in a harmonic analysis context (see the book \cite{Meyer}), which he termed harmonious sets, and observed that some specific examples of Meyer sets that he called \textit{model sets} could be produced through an elementary geometric method. More than a decade latter and independently of Meyer's work, the geometric method to produce model sets have been reintroduced by several authors in an attempt to answers a problem in solid state physics, by producing point patterns of $\mathbb{R}^{2}$ and $\mathbb{R}^3$ having a diffraction picture displaying rotational symmetries that cannot hold for fully periodic point sets, as it was first realized by considering the vertex sets of the renowned Penrose tilings of the plane \cite{Mac}. Nowadays many interesting examples of model sets and tilings affiliated with model sets can be found in the literature \cite{BaakeGrimm}, illustrating the elegance and fruitfulness of Meyer's method.\\


It has early been observed that a certain property holds for model sets when one looks at the dynamical system canonically associated to any point set: It is always an \textit{almost 1-to-1 cover} of some Kronecker system. Such property was first recognized by Robinson \cite{Rob}, who asserted that the dynamical system of a Penrose tiling is a minimal almost 1--1 extension of an $\mathbb{R}^2$-action over the torus $\mathbb{T}^4$ by rotation. This may be explained by the existence of an underlying model set structure for such tilings, and by a direct use of the model set construction such this has later been obtained for the Fibonacci model set of the line \cite{BaaHerPle}, whose dynamical system was shown to be a minimal almost 1--1 extension of an $\mathbb{R}$-action by rotation over the torus $\mathbb{T}^2$. These results were generalized for certain repetitive model sets of a Euclidean space by Forrest, Hunton and Kellendonk \cite{ForHunKel}, and independently by Schlottmann \cite{Sch2} for arbitrary repetitive model sets of a $\sigma$-compact, locally compact Abelian group, where in both works it is shown that the model set construction immediately describes the dynamical system of a repetitive model set as a minimal almost 1--1 extension of a Kronecker system (this is the so-called \textit{torus parametrization} of a model set). Dynamical systems with this property appeared long before in general topological dynamics, and in particular in the work of Veech \cite{Vee} who termed them \textit{almost automorphic}.\\


It was thus reasonable to think that this property should characterize on a dynamical side (repetitive) model sets among the wider class of (repetitive) Meyer sets, and a proof of such characterization came with the work of by Baake, Lenz and Moody \cite{BaaLenMoo} and Lee and Moody \cite{LeeMoo}, who proved using the geometric construction developed in \cite{BaaMoo0} that repetitive model sets displaying the additional condition of \textit{regularity} are characterized by a strong form of almost automorphy on their associated dynamical system. Although regularity is a very reasonable condition notably concerning the diffraction theory of point sets \cite{Hof,Sch2}, it is not necessary for the almost automorphy condition to appear, and since regularity has a crucial task in the arguments of \cite{BaaLenMoo,LeeMoo} a dynamical characterization that applies for model sets not necessarily regular remained an interesting open problem. This has been obtained using different arguments in \cite{Auj2}, where it is proved in the case $\mathbb{G}= \mathbb{R}^d$ the following:\\

\textbf{Theorem.} \textit{A repetitive model set is precisely a repetitive Meyer set whose associated dynamical system is almost automorphic.}\\

The proof set in \cite{Auj2} involves Ellis semigroup and regional proximality for Meyer sets, two objects coming from abstract topological dynamics \cite{Auslander}, for which a detailed account in the context of point sets can be found in \cite{AujBarKelLen}. It is also detailed in \cite{AujBarKelLen} how  the result of \cite{BaaLenMoo,LeeMoo} and the above theorem are part of a natural hierarchical description of repetitive Meyer sets involving their associated dynamical systems.

In this note, we aim to revisit this result by providing a greatly simplified proof, which is shown here to hold for point sets in arbitrary $2^{nd}$ countable, $\sigma$-compact locally compact Abelian groups. Our starting point here is a result of schlottmann upon characterizing repetitive model sets by an intrinsic condition on the return times of the set \cite{Sch}, that is unexpectedly very close to almost automorphy. In particular we do not make use of the torus parametrization anywhere in this work. The proof we give here consists of a combination of Schlottmann's original result (Theorem \ref{theo:Schlottmann}) with a general statement on almost automorphy (Theorem \ref{theo:Veech}), and with an intermediate statement connecting the two (Proposition \ref{prop:additivity}).

\section{Meyer sets}\label{section.Meyer.sets}

In this note we consider Meyer sets of a locally compact Abelian (LCA) group $\mathbb{G}$, which for the sake of simplicity is assumed  $\sigma$-compact and $2^{st}$ countable. One of the various definitions of a Meyer set is the following: A Meyer set $\Lambda$ of $\mathbb{G}$ is a subset that is both uniformly discrete and relatively dense in $\mathbb{G}$, and for which there exists a finite set $F$ such that $ \Lambda - \Lambda \, \subseteq \, \Lambda +F$. This definition admits many different, though equivalent, formulations, which the reader can find in \cite{Moo} in the case $\mathbb{G} = \mathbb{R}^d$ and \cite{Str3} in our slightly more general setting.\\

A trivial example of Meyer set is a uniformly discrete and relatively dense subgroup $\Gamma $, a lattice, of $ \mathbb{G}$, and any sum $\Gamma + F$ with $F$ finite.	Much less trivial examples of Meyer sets are the \textit{model sets}, also called Cut $\&$ Project sets: To construct a model set one begins by considering a product $H \times \mathbb{G}$ with a LCA group $H$ together with a lattice $\Gamma \subset H \times \mathbb{G}$ in it. The triple $(H, \Gamma ,  \mathbb{G})$ is called a cut $\&$ project scheme. Then choosing a compact subset $W \subseteq H$ that is supposed equal to the closure of its interior one ends up with the point set
\[ \mathfrak{P}(W) := \left\lbrace \, t \in \mathbb{G} \, : \, \exists \, (t',t) \in \Gamma , \, t' \in W \,\right\rbrace 
\]

\vspace{0.1cm}
One usually suppose that $\Gamma$ projects in a 1--1 way into $\mathbb{G}$, but this can be achieved by modding out the discrete subgroup $\Gamma _0 :=\Gamma \cap H \times \left\lbrace 0\right\rbrace $ in $H$ and consider the cut and project scheme $(H/_{\Gamma _0}, \Gamma /_{\Gamma _0}, \mathbb{G})$. It is moreover often required that the projection of $\Gamma$ in $H$ has dense image, which is obtained by taking the closure $H'$ of this range in replacement of the whole $H$. None of these operations alters the model sets as given above.\\

There are natural generalizations of this definition: For instance, instead of requiring $W$ to have dense interior we could ask this interior to be only non-empty, or even empty at all \cite{Moo3}. One also may formulate this in terms of the indicator function of the subset $W$ and then replace it by different types of function (continuous \cite{LenRic,Str3}, Riemann measurable, Baire class 1). Ultimately, one may also drop the commutativity requirement on all involved groups. In our case, we will be particularly interested in a certain type of model sets: A model set $\mathfrak{P}(W)$ is called \textit{non-singular} when there is no point of $\Gamma$ whose projection in $H$ belongs to the boundary of $W$. In other terms, a model set is non-singular when it satisfies
\vspace{0.1cm}
\[ \mathfrak{P}(\mathring{W})  = \mathfrak{P}(W)
\]

\vspace{0.1cm}

\section{The dynamical system of a Meyer set}\label{section:Meyer.set.dynamical.system}
We wish to characterize the Meyer sets that are non-singular model sets in terms of a certain dynamical system canonically associated to each of them. Given a Meyer set $\la$, the construction of its dynamical system \xlag goes as follows: We consider the set $X_\la$, the hull of $\la$, of all Meyer sets in which any finite subset appears somewhere in $\la$,
\[ X_{\mathsf{\Lambda}}:= \left\lbrace \, \Lambda  \subset \mathbb{G} \, : \,\forall \,K\subset \mathbb{G} \text{ compact } \, \exists \, t\in \mathbb{G},\, \Lambda  \cap K= (\la-t)\cap K \,  \right\rbrace 
\]
This set is endowed with a topology that comes from a uniformity (see \cite[Chapter 6]{Kelley} for background), for which a basis is given, for an arbitrary compact set $K$ of $\mathbb{G}$ and $U$ any neighborhood of $0$ in $\mathbb{G}$, by the sets
\[ \mathcal{U}_{K,U} := \left\lbrace \, (\Lambda , \Lambda ')\in X _\la \times X _\la \, : \, \exists \, t, t' \in U, \,  (\Lambda -t) \cap K \equiv (\Lambda '-t') \cap K\right\rbrace 
\]
The set $X _\la$ endowed with the topology coming from this uniformity is a compact Hausdorff space \cite{Sch2}, and since $\mathbb{G}$ is assumed $\sigma$-compact and $2^{nd}$ countable this uniformity has a countable basis and thus is metrizable \cite[Chapter 6, Theorem 13]{Kelley}. Finally, it is naturally endowed with a $\mathbb{G}$-action by $t.\Lambda := \Lambda -t$ and this action is by homeomorphism. By construction, the Meyer set $\la$ now viewed as an element of $X_\la$ is a transitive element, that is, its $\mathbb{G}$-orbit is dense in the dynamical system $(X_{\mathsf{\Lambda}}, \mathbb{G})$. \\

 In this note we will exclusively be interested in Meyer sets that are repetitive: A Meyer set $\la$ is repetitive whenever for any compact subset $K$ of $\mathbb{G}$ the following set
\[ \mathsf{\Lambda}_ K := \left\lbrace \, t\in \mathbb{G} \;  : \; \displaystyle (\mathsf{\Lambda} -t)\cap K = \mathsf{\Lambda}\cap K \right\rbrace 
\]
is relatively dense in $\mathbb{G}$. It admits an elegant interpretation in terms of the associated dynamical system \xlag as follows. Observe first that for a compatible metric $d$ on $X_\la$, for any compact set $K$ and open neighborhood $U$ of $0$ in $\mathbb{G}$ there must exists an $\varepsilon > 0 $ such that whenever $d(\Lambda, \Lambda ') < \varepsilon$ then $(\Lambda , \Lambda ') \in \mathcal{U}_{K, U}$, and conversely given $\varepsilon > 0$ there exists a pair of a compact set $K$ and an open neighborhood $U$ of $0$ in $\mathbb{G}$ such that the reverse implication holds. Now let us introduce the collection of sets given for $\varepsilon > 0$ by
\[ P_\varepsilon := \left\lbrace \, t\in \mathbb{G} \;  : \;  d(t.\la , \la ) < \varepsilon \right\rbrace 
\]
What we can observe is then that for any pair of a compact set $K$ and an open neighborhood $U$ of $0$ in $\mathbb{G}$ there exists a $\varepsilon > 0$ such that $P_\varepsilon \subseteq \mathsf{\Lambda}_ K +U$, and conversely, for any $\varepsilon > 0$ there are $K$ and $U$ as before such that $ \mathsf{\Lambda}_ K +U\subseteq P_\varepsilon $. From this observation we conclude that the Meyer set $\la$ is repetitive precisely when the $P_\varepsilon$'s are relatively dense in $\mathbb{G}$, that is, when $\la$ is an almost periodic point in the dynamical system $(X_{\mathsf{\Lambda}}, \mathbb{G})$. It is a general fact of topological dynamics that this happens if and only if the dynamical system \xlag is minimal.\\ 

Thus, when a Meyer set $\la$ is repetitive its dynamical system \xlag is equally generated by any other Meyer in it, in which case we may hide the subscript $\la$ and simply write \xg, and refer it as a \textit{minimal system of Meyer sets}. When a minimal system of Meyer sets \xg contains a non-singular model set, or equivalently is generated by a non-singular model set (these are always repetitive) then we call it a \textit{minimal system of model sets}, and call its elements \textit{repetitive model sets}. For more details see for instance \cite{Auj2}.

\vspace{0.1cm}
\section{Schlottmann's characterization of model sets}\label{section.Schlottmann}

Let us recall the result of Schlottmann:

\vspace{0.3cm}
\begin{theorem}\label{theo:Schlottmann}\cite[Theorem 2]{Sch} For a Meyer set $\la $ the following statements are equivalent:

\vspace{0.2cm}
\begin{itemize}
\item[1.] $\la$ is a non-singular model set.

\vspace{0.3cm}
\item[2.] Each $\mathsf{\Lambda}_ K$ is relatively dense, and for any compact $K$ there is a compact $K'$ with $$ \mathsf{\Lambda}_ {K'} - \mathsf{\Lambda}_ {K'} \subseteq \mathsf{\Lambda}_ K $$

\end{itemize}
\end{theorem}

\vspace{0.1cm}

The complete proof of this can be found in Schlottmann's original paper \cite{Sch}. We warn the reader that non-singularity is termed "regularity" in \cite{Sch}, which is unrelated with nowadays' notion of regularity used for model sets. We shall see now that, given a Meyer set $\la$ with associated dynamical system \xlag ,  Schlottmann's characterization can be expressed using the collection $(P_\varepsilon ) _\varepsilon $ introduced in Section \ref{section:Meyer.set.dynamical.system} for some compatible metric $d$ on $X_\la$: 

\vspace{0.2cm}
\begin{proposition}\label{prop:additivity} For $\la$ a Meyer set and \xlag its associated dynamical system the following two statements are equivalent:

\vspace{0.1cm}
\begin{itemize}
\item[1.] Each $\mathsf{\Lambda}_ K$ is relatively dense, and for a compact $K$ there is a compact $K'$ with $$ \mathsf{\Lambda}_ {K'} - \mathsf{\Lambda}_ {K'} \subseteq \mathsf{\Lambda}_ K $$

\item[2.] Each $P_{\varepsilon }$ is relatively dense, and for each $\varepsilon > 0$ there exists some $\delta >0$ such that $$  P_{\delta } - P_{\delta } \subseteq P_{\varepsilon }$$
\end{itemize}
\end{proposition}

\vspace{0.1cm}

\begin{proof} We already discussed the equivalence between relative denseness of the $\la _K$'s and of the $P_\varepsilon$' in Section \ref{section:Meyer.set.dynamical.system}. Now one easily shows that point 1 implies point 2: Since $d$ is compatible, from the discussion of Section \ref{section:Meyer.set.dynamical.system} we have for any given $\varepsilon > 0$ a compact $K$ and and a neighborhood $U$ of $0$ in $\mathbb{G}$ with $ \mathsf{\Lambda}_ K +U\subseteq P_\varepsilon $, and taking $K'$ and $U'$ such that $ \mathsf{\Lambda}_ {K'} - \mathsf{\Lambda}_ {K'} \subseteq \mathsf{\Lambda}_ K $ and $U'-U' \subseteq U$, and then choosing $\delta > 0$ with $P_\delta \subseteq \mathsf{\Lambda}_{K'} +U'$ gives $P_{\delta } - P_{\delta } \subseteq P_{\varepsilon }$, as desired.\\

Let us show that conversely point 2 implies point 1: Again from the discussion of Section \ref{section:Meyer.set.dynamical.system} we deduce that for $K$ chosen and any open neighborhood $U$ of $0\in \mathbb{G}$ there exists, using our assumption, a $K'$ compact as well as another open neighborhood $U'$ of $0\in \mathbb{G}$, that can be chosen arbitrarily small, such that $(\mathsf{\Lambda}_ {K'} +U') - (\mathsf{\Lambda}_ {K'} +U') \subseteq \mathsf{\Lambda}_ {K} +U$. This yields in particular $\mathsf{\Lambda}_ {K'} - \mathsf{\Lambda}_ {K'} \subseteq \mathsf{\Lambda}_ {K} +U$. Now let us observe that for a compact $K$ that intersect $\la$ the point set $\la _{K}$ belongs to the difference set $\la - \la$: Indeed having $t \in \la _{K}$ means that $\la \cap K$ and $(\la -t)\cap K$ are equal, and selecting some point $p$ in it ensures that both $p$ and $p+t$ belong to $\la$, providing $t = (p+t) - p \in \la - \la$. Thus, applying this to $K'$ yields that $\mathsf{\Lambda}_ {K'} - \mathsf{\Lambda}_ {K'} $ belongs to $\la - \la + \la - \la$ independently on $K'$, provided it is large enough to intersect $\la$. Now as we can freely reduce $U$ (although this would change $K'$) we may select it small enough so that $(\la - \la + \la - \la ) \cap (\mathsf{\Lambda}_ {K} +U)$ is actually contained in $\mathsf{\Lambda}_ {K}$: Indeed one can check that any  open set $U$ whose intersection with $\la - \la + \la - \la + \la - \la$ is $ \left\lbrace 0\right\rbrace $ works, and such an open set always exists from the Meyer property of $\mathsf{\Lambda}$. Thus, for such a choice of open $U$ we end up with a compact $K'$ that satisfies the desired inclusion $\mathsf{\Lambda}_ {K'} - \mathsf{\Lambda}_ {K'} \subseteq \mathsf{\Lambda}_ {K}$.
\end{proof}

\vspace{0.2cm}

\section{A general result on almost automorphy}\label{Section.Veech}

In this section we consider an abstract dynamical system $(X, \mathbb{G})$, where $X$ is a compact metric space and $\mathbb{G}$ a LCA group that acts jointly continuously on $X$. We moreover assume the existence of a transitive point $\mathrm{x}\in X$, that is, a point whose $\mathbb{G}$-orbit is dense in $X$.\\

 Recall that a factor map $\pi :X \longrightarrow Y$ from a dynamical system \xg to another dynamical system \yg is a continuous onto mapping that is a $\mathbb{G}$-map, that is, such that $\pi (t.x )  = t.\pi (x)$ occurs for any $x\in X$ and any $t\in \mathbb{G}$. A dynamical system \xg always admit a greatest factor whose $\mathbb{G}$-action is equicontinuous, its so-called \textit{maximal equicontinuous factor} $(X_{eq}, \mathbb{G})$, which is uniquely defined up to conjugacy, see \cite{AujBarKelLen} for details. The factor map from a dynamical system onto its maximal equicontinuous factor is in particular unique up to left-composing the factor map with a $\mathbb{G}$-commuting homeomorphism of $X_{eq}$.
 
 For a dynamical system \xg with a transitive point $\mathrm{x}$ the maximal equicontinuous factor is a minimal \textit{Kronecker system}: The space $X_{eq}$ admits a compact (Hausdorff) Abelian group structure, and there exists a continuous group morphism $r : \mathbb{G} \longrightarrow X_{eq}$ with dense range such that the $\mathbb{G}$-action on $X_{eq}$ is given by addition $t. z := r(t) + z$ for any $z\in X_{eq}$ and $t\in \mathbb{G}$. The group structure on $X_{eq}$ is not canonical here. In fact, any element $z\in X_{eq}$ can play the role of unit, and this is with respect to a unique compact Abelian group structure on $X_{eq}$. Once an element $z \in X_{eq}$ is chosen as unit, the morphism $r$ is recovered through the $\mathbb{G}$-action on $X_{eq}$ by $r(t) := t.z$ for any $t\in \mathbb{G}$.\\
 
A dynamical system \xg is \textit{almost automorphic} when its factor map $\pi : X \longrightarrow X_{eq}$ is 1--1 over some element of $X_{eq}$. A point $x\in X$ where $\pi $ is 1--1 is called an almost automorphic point ("non-singular point" would have also been a good terminology). In this section, we shall characterize the property of almost automorphy of a dynamical system \xg and of a point $\mathrm{x}\in X$ by means of the sets, for any $\varepsilon > 0$, of $\varepsilon$-return times of $\mathrm{x}$
\[ P_\varepsilon := \left\lbrace \, t\in \mathbb{G} \;  : \;  d(t.\mathrm{x}, \mathrm{x}) < \varepsilon \right\rbrace \]
Each of these sets contains the origin $0 \in \mathbb{G}$. The following result shows that almost automorphy of the system \xg corresponds to a simple algebraic structure on the set of return times of a single transitive point:

\vspace{0.3cm}
\begin{theorem}\label{theo:Veech} Let \xg be a dynamical system and $\mathrm{x} \in X$ a transitive point. Then the following statements are equivalent:

\vspace{0.2cm}
\begin{itemize}
\item[1.] Each $P_{\varepsilon }$ is relatively dense, and for each $\varepsilon > 0$ there exists some $\delta >0$ such that $$  P_{\delta } - P_{\delta } \subseteq P_{\varepsilon }$$

\item[2.] \xg is almost automorphic, with $\mathrm{x} \in X$ an almost automorphic point.

\end{itemize}
\end{theorem}

\vspace{0.3cm}


\begin{proof}
Assume that point 1 holds: The collection of subsets $(P_\varepsilon )_{\varepsilon > 0}$ of $\mathbb{G}$ defines the basis $(B_\varepsilon )_{\varepsilon > 0 }$ of a translation-invariant and inversion-invariant uniformity on $\mathbb{G}$, where $B_\varepsilon \subseteq \mathbb{G} \times \mathbb{G}$ is the set of pairs $(t, t')$ with $t-t' \in P_\varepsilon$. The topology coming from this uniformity, which we call here the $(P_\varepsilon)_\varepsilon $-topology, turns $\mathbb{G}$ into an Abelian (possibly non Hausdorff) topological group, for which $(P_\varepsilon)_\varepsilon $ is a neighborhood basis of $0 \in \mathbb{G}$. This topological group admits a completion by an Abelian (Hausdorff) compete topological group $\mathcal{Z}$ with completion morphism $r : \mathbb{G} \longrightarrow \mathcal{Z}$, and
\begin{itemize}

\vspace{0.2cm}
\item[•] with dense range in $\mathcal{Z}$,
\item[•] having kernel $Ker (r) = \bigcap _{\varepsilon > 0 } P_\varepsilon $,
\item[•] which is continuous and open for the $(P_\varepsilon)_\varepsilon $-topology on $\mathbb{G}$
\end{itemize}

\vspace{0.2cm}
For details one may see \cite[Appendix A]{Sch}. Then the mapping defined along the $\mathbb{G}$-orbit of $ \mathrm{x} \in X$ by $t.\mathrm{x} \longmapsto r(t)\in \mathcal{Z}$ is well-defined, because each time we get $t. \mathrm{x} = t'. \mathrm{x}$ then $t-t'$ belongs to any $P_\varepsilon$ and thus to the kernel of the morphism $r$, giving $r(t) = r(t')$. It is moreover uniformly continuous and therefore extends in a continuous map $\pi : X \longrightarrow \mathcal{Z}$: For, we need to show that for any $\varepsilon > 0$ there exists a $\delta > 0$ such that whenever $d(t. \mathrm{x}, t'. \mathrm{x} ) < \delta $ then $t-t' \in P_\varepsilon$. To see this, observe that for $\varepsilon > 0$ there exists by assumption a $\varepsilon ' > 0$ such that $ P_{\varepsilon ' } - P_{\varepsilon '} \subseteq P_{\varepsilon }$. Now because the $(P_\varepsilon ) _\varepsilon $ are relatively dense the transitive point $\mathrm{x}$ is almost periodic, i.e. \xg is minimal, so the set of translates $s. B (\mathrm{x}, \varepsilon ')$, $s\in \mathbb{G}$, covers $X$. One is thus able to extract a finite subcover given by $s\in F$ with $F$ finite, which admits a Lebesgue number $\delta > 0$. This guarantees that if one has $d(t. \mathrm{x}, t'. \mathrm{x} ) < \delta $ then there is an $s\in F$ such that $t. \mathrm{x}$ and $  t'.\mathrm{x} $ belong to $s.B(\mathrm{x}, \varepsilon ')$, and as a result both $(t-s).\mathrm{x}$ and $(t'-s).\mathrm{x}$ belong to the ball $B(\mathrm{x}, \varepsilon ' )$. Hence $t-t' = (t-s) -(t'-s)\in P_{\varepsilon '} - P_{\varepsilon '} \subseteq P _\varepsilon$, as desired. The resulting map $\pi$ satisfies $\pi (\mathrm{x}) = \mathfrak{o}\in \mathcal{Z}$, the unit element in $\mathcal{Z}$, and is onto since its image is compact and contains the range of $r$. It is naturally a $\mathbb{G}$-map as one may check that $\pi (t. x) = r(t). \pi (x)$ holds for any $x\in X$ and $t\in \mathbb{G}$. Therefore, $\pi$ is a factor map form \xg to the Kronecker system $(\mathcal{Z}, \mathbb{G})$ with Kronecker action given by the morphism $r$, which is continuous since is given by $r(t) = \pi (t.\mathrm{x})$ for each $t\in \mathbb{G}$. It now suffices to show that $ \mathrm{x} $ is a one-point fiber with respect to this factor map to conclude that $(\mathcal{Z}, \mathbb{G})$ must be the maximal equicontinuous factor of \xg and that point 2 holds. Let us show this: Suppose that $x\in X$ satisfies $\pi (x) = \pi (\mathrm{x})$. Since $\mathrm{x}$ is transitive there is a sequence $(t_n)_n \subset \mathbb{G}$ such that $t_n.\mathrm{x}$ converges to $x$, and this gives in the factor $\mathcal{Z}$ the convergence of $r(t_n)$ towards the unit $\mathfrak{o}\in \mathcal{Z}$. Since $r : \mathbb{G} \longrightarrow \mathcal{Z}$ is open when $\mathbb{G}$ is equipped with the $(P_\varepsilon)_\varepsilon $-topology and the collection $(P_\varepsilon)_\varepsilon $ forms a neighborhood basis of $0 \in \mathbb{G}$, we deduce that for any $\varepsilon > 0$ the sequence $(t_n)$ must belong to $P_\varepsilon $ for $n$ great enough. But by construction of the sets $P_\varepsilon$ this means that $d(t_n.\mathrm{x}, \mathrm{x}) < \varepsilon $ for any $\varepsilon > 0$ and $n$ great enough. This latter argument guarantees that $x = \mathrm{x}$, giving point 2.\\

Assume that point 2 holds: The dynamical system \xg admits a factor map $\pi : X \longrightarrow X_{eq} $ that is 1--1 at $\mathrm{x}$, where $(X_{eq}, \mathbb{G})$ is the maximal equicontinuous factor with Kronecker action given by a continuous morphism $r : \mathbb{G} \longrightarrow X_{eq}$. Since $\pi$ is 1--1 at $\mathrm{x}$ one easily shows (by contradiction) that for any $\varepsilon > 0$ there is a neighborhood $U_\varepsilon$ of $\pi (\mathrm{x}) \in X_{eq}$ such that $\pi ^{-1}(U_\varepsilon) \subseteq B_d(\mathrm{x}, \varepsilon )$. From this we deduce that the sets $P_\varepsilon $ are relatively dense in $\mathbb{G}$ (e.g. $\mathrm{x}$ is an almost periodic point): Indeed the point $\pi (\mathrm{x})$ is transitive and thus almost periodic in the the Kronecker system $(X_{eq}, \mathbb{G})$, so the set of $t\in \mathbb{G}$ such that $r(t) +\pi (\mathrm{x})$ belongs to $U_\varepsilon$ is relatively dense in $\mathbb{G}$. Each such point admits $t. \mathrm{x}$ as unique lift in $X$, which thus belong to $B_d(\mathrm{x}, \varepsilon )$ from the very choice of $U_\varepsilon$, giving that $t\in P_\varepsilon$, so this latter is relatively dense. Now since $X_{eq}$ is a topological group one may choose for each $\varepsilon > 0$ a neighborhood $ U'_\varepsilon$ of the unit $\mathfrak{0}\in X_{eq}$ such that $U'_\varepsilon -U'_\varepsilon \subseteq U_\varepsilon - \pi (\mathrm{x})$. As $\pi$ is continuous, there is some $\delta >0$ such that $\pi (B_d(\mathrm{x}, \delta  )) \subseteq \pi (\mathrm{x}) +U'_\varepsilon$, and therefore whenever $t, t'\in P_{\delta }$ then $t.\mathrm{x}$ and $t'.\mathrm{x}$ have images $r(t) +\pi (\mathrm{x})$ and $r(t') +\pi (\mathrm{x})$ that belong to $\pi (\mathrm{x}) +U'_\varepsilon$. Thus $r(t-t')$ lies in $U'_\varepsilon-U'_\varepsilon \subseteq U_\varepsilon -\pi (\mathrm{x})$, and therefore $ r(t-t') +\pi (\mathrm{x})$ belongs to $ U_\varepsilon$. Consequently its lift $(t-t').\mathrm{x}  $ lies in $B_d(\mathrm{x}, \varepsilon )$, yielding $t-t' \in P_\varepsilon$. This gives point 1.
\end{proof}

\vspace{0.3cm}
\begin{remark} It is possible to set an elementary proof of Theorem \ref{theo:Veech} by invoking a much more general (and deeper) result due to Veech \cite{Vee2}, which may be stated as follows. First, since both statements in Theorem \ref{theo:Veech} imply minimality of the system \xg we can assume this condition. Then consider a metric minimal system $(X,\mathbb{G})$, with factor map $\pi : X \longrightarrow X_{eq} $ over its maximal equicontinuous factor, and a point $\mathrm{x} \in X$ with associated family of sets $(P_\varepsilon )_{\varepsilon > 0}$. The result of \cite{Vee2} then describes the set of points in $X$ belonging in the fiber of $\mathrm{x}$ over $X_{eq}$, as 
\[ \pi^{-1}(\pi(\mathrm{x}))= \bigcap _{\varepsilon > 0} \overline{(P_\varepsilon - P_\varepsilon).\mathrm{x}} \]
Using this result one easily shows Theorem \ref{theo:Veech}: Indeed, if point 1 holds then the fiber $\pi^{-1}(\pi(\mathrm{x}))$ is equally given by the intersection of the sets $\overline{P_\varepsilon.\mathrm{x}}$, which is nothing but the singleton $\left\lbrace \mathrm{x}\right\rbrace $, so that point 2 holds. Conversely if point 2 holds then for any $\varepsilon > 0$ there is a $\delta > 0$ such that $(P_\delta - P_\delta).\mathrm{x}$ must belong to $B(x, \varepsilon )$, giving $P_\delta - P_\delta \subseteq P_\varepsilon$ and thus point 1 holds.

 As a result, the characterization of model sets by dynamical systems as presented here is mainly the combination of Schlottmann's characterization and of Veech's result on the equicontinuous structure relation for minimal systems. Despite the existence of this very direct argument, we decided to keep the proof of Theorem \ref{theo:Veech} as presented above, first for self-containment of the note, and second because it involves a completion process of topological groups that is very similar to the one appearing in the proof of Schlottmann's characterization given in \cite{Sch}.
\end{remark}

\vspace{0.4cm}

\textbf{Acknowledgments.}  It is a pleasure to thank Daniel Lenz for encouraging me to write this note, as well as Johannes Kellendonk and Christoph Richard for their valuable comments and for pointing out several references.

\vspace{0.3cm}

\bibliographystyle{plain}

\bibliography{Biblio}

\end{document}